\documentclass[12pt,reqno]{amsart}
\usepackage{fullpage,a4wide}

\usepackage{graphicx}
\usepackage[draft]{hyperref}
\usepackage{amsmath,amsopn,amssymb,amsfonts,stmaryrd}
\usepackage{verbatim}
\usepackage{amsthm}
\usepackage{mathtools}
\usepackage{color}
\usepackage{enumitem}
\usepackage[framemethod=TikZ]{mdframed}
\usepackage{bbm}
\usepackage{mathrsfs}
\usepackage{booktabs}
\usepackage{caption}
\usepackage{bm}
\usepackage{tensor}
\usepackage{cleveref}

\renewcommand{\S}{\mathbb{S}}

\newcommand{\CC}{\mathcal{C}}

\newcommand{\CE}{\mathcal{E}}

\newcommand{\CS}{\mathcal{S}}

\newcommand{\N}{\mathbb N}
\newcommand{\C}{\mathbb C}

\newcommand{\ssum}{{\textstyle{\displaystyle\sum}}}

\theoremstyle{plain}
\newtheorem{thm}{Theorem}[section]

\newtheorem{lem}[thm]{Lemma}

\newtheorem{conj}[thm]{Conjecture}

\theoremstyle{definition}

\numberwithin{equation}{section}

\newcommand{\Flo}[1]{\left\lfloor #1\right\rfloor}

\newcommand{\pa}[2]{\left(\frac{#1}{#2}\right)}
\newcommand{\tpa}[2]{\left(\tfrac{#1}{#2}\right)}

\def\d{\delta}

\def\l{\lambda}

\def\w{\omega}

\def\z{\zeta}

\def\vth{\vartheta}

\def\e{\varepsilon}

\def\d{\delta}

\def\l{\lambda}

\def\w{\omega}

\def\z{\zeta}

\def\vth{\vartheta}

\def\e{\varepsilon}

\newcommand{\re}{{\rm Re}}

\newcommand{\err}{\mathrm{err}}

\newcommand{\R}{\mathbb R}
\newcommand{\Z}{\mathbb Z}

\newcommand{\lrb}[1]{\left(#1\right)}

\newcommand{\rank}{\mathrm{rank}}
\newcommand{\crank}{\mathrm{crank}}

\setlist[itemize]{noitemsep, topsep=0pt}

\allowdisplaybreaks

\makeatletter
\newcommand{\vast}{\bBigg@{2}}
\newcommand{\Vast}{\bBigg@{5}}
\makeatother

\newcommand{\Pmod}[1]{\ \, ( \mathrm{mod} \, #1 )}

\makeatletter
\@namedef{subjclassname@2020}{%
	\textup{2020} Mathematics Subject Classification}
\makeatother

\title{Biases among classes of rank-crank partitions$\bm{\Pmod{11}}$}
\author{Kathrin Bringmann}
\author{Badri Vishal Pandey}
\address{University of Cologne, Department of Mathematics and Computer Science, Weyertal 86-90, 50931 Cologne, Germany}
\email{kbringma@math.uni-koeln.de}
\email{bpandey@uni-koeln.de, badrivishal9451@gmail.com}


\begin{document}
\maketitle

\begin{abstract}
	In this paper, we prove inequalities for ranks, cranks, and partitions among different classes modulo 11. These were conjectured by Borozenets.
\end{abstract}

\section{Introduction and statement of results}

We start by recalling Ramanujan's famous congruences for the partition function $p(n)$. He proved in \cite{R} that for $n\in\N_0$,
\begin{equation}\label{E:Dysoncong}
	p(5n+4) \equiv 0\Pmod5,\quad p(7n+5) \equiv 0\Pmod7,\quad p(11n+6) \equiv 0\Pmod{11}.
\end{equation}
To combinatorically explain the first two congruences, Dyson \cite{D} introduced the {\it rank} of a partition $\l$. This is given by
\[
	\rank(\l) := \text{largest part of }\l - \text{number of part of }\l.
\]
Dyson also conjectured the existence of another partition statistic that would explain all of the three congruences, which he called ``crank". Garvan discovered the so-called vector crank in \cite{G}, which was subsequently reformulated by Andrews and Garvan as the partition statistic crank \cite{AG}. This successfully completed the search for combinatorial decompositions of the three congruences in \eqref{E:Dysoncong}. If $o(\l)$ denotes the number of ones in $\l$, and $\mu(\l)$ is the number of parts strictly larger than $o(\l)$, then the {\it crank} is defined as
\[
	\crank(\l) :=
	\begin{cases}
		\text{largest part of }\l & \text{if }o(\l)=0,\\
		\mu(\l)-o(\l) & \text{if }o(\l)>0.
	\end{cases}
\]

Let $M(m,n)$ (resp. $N(m,n)$) denote the number of partitions of $n$ with crank (resp. rank) $m$. Moreover denote by $M(a,c;n)$ (resp. $N(a,c;n)$) the number of partitions of $n$ with crank (resp. rank) congruent to $a\Pmod c$. There are many inequalities known about ranks and cranks. For example Andrews and Lewis \cite{AL} proved several inequalities and conjectured the following:
\begin{align*}
	M(0,3;3n) &> M(1,3;3n),\qquad M(0,3;3n+1) < M(1,3;3n+1),\\
	M(0,3;3n+2) &< M(1,3;3n+2) \quad\text{unless } n\in\{1,4,5\},\\
	N(0,3;n) &< N(1,3;n) \quad\text{if } n\equiv0,2\Pmod3,\\
	N(0,3;n) &> N(1,3,n) \quad\text{if } n\equiv1\Pmod3.
\end{align*}
The inequalities for the crank have been shown by Kane \cite{K} and the ones for the rank by the first author \cite{Br}. Since then there has been many studies of such inequalities (see e.g. \cite{BK,CM}).

Recently, Borozenets \cite{Bo} considered rank and crank statistic modulo $11$, where he proved many inequalities among different classes of rank and crank modulo $11$ as well as conjectured some new inequalities among them. The goal of this paper is to prove the following conjectures made in \cite{Bo}. To state them, throughout, for sequences $a=\{a(n)\}_{n=1}^\infty$ and $b=\{b(n)\}_{n=1}^\infty$, we use the notation $a\le_N b$ (resp. $a\le b$) to mean that $a(n)\le b(n)$ for $n\ge N$ (resp. $n\geq 0$). For $0\le j,d<11$, define $N_{j,d}(n):=N(j,11;11n+d)$, $M_{j,d}(n):=M(j,11;11n+d)$, and $p_d(n):=p(11n+d)$. 

\begin{conj}[Conjecture 6.15 of \cite{Bo}]\label{CONJ:rankcrankconj}
	We have the following inequalities
	\begin{align*}
		N_{5,0} \le N_{4,0} \le_2 N_{3,0} \le M_{1,0} &\le \frac{p_0}{11} \le M_{0,0} \le_1 N_{2,0} \le N_{1,0} \le_3 N_{0,0},\\
		N_{5,1} \le N_{4,1} \le N_{3,1} \le_1 M_{0,1} \le M_{2,1} &\le \frac{p_1}{11} \le M_{1,1} \le_1 N_{2,1} \le N_{1,1} \le N_{0,1},\\
		N_{5,2} \le N_{4,2} \le N_{3,2} \le M_{0,2} &\le \frac{p_2}{11} \le M_{2,2} \le_1 N_{2,2} \le N_{1,2} \le_3 N_{0,2},\\
		N_{5,3} \le N_{4,3} \le N_{3,3} \le M_{1,3} &\le \frac{p_3}{11} \le M_{0,3} \le N_{2,3} \le_1 N_{1,3} \le_2 N_{0,3},\\
		N_{5,4} \le_{3} N_{4,4} \le N_{3,4} \le_1 M_{1,4} &\le \frac{p_4}{11} \le M_{0,4} \le_1 N_{2,4} \le N_{1,4} \le_3 N_{0,4},\\
		N_{5,5} \le N_{4,5} \le_1 N_{3,5} \le M_{2,5} &\le \frac{p_5}{11} \le M_{0,5} \le N_{2,5} \le N_{1,5} \le N_{0,5},\\
		N_{5,6} \le_1 N_{4,6} \le N_{3,6} &\le \frac{p_6}{11} \le N_{2,6} \le N_{1,6} \le_1 N_{0,6},\\
		N_{5,7} \le N_{4,7} \le N_{3,7} \le M_{0,7} &\le \frac{p_7}{11} \le M_{1,7} \le N_{2,7} \le_1 N_{1,7} \le N_{0,7},\\
		N_{5,8} \le N_{4,8} \le N_{3,8} \le M_{0,8} &\le_3 \frac{p_8}{11} \le_3 M_{1,8} \le N_{2,8} \le N_{1,8} \le_3 N_{0,8},\\
		N_{5,9} \le N_{4,9} \le N_{3,9} \le_1 M_{0,9} &\le \frac{p_9}{11} \le M_{1,9} \le N_{2,9} \le N_{1,9} \le_2 N_{0,9},\\
		N_{5,10} \le N_{4,10} \le N_{3,10} \le_1 M_{3,10} &\le \frac{p_{10}}{11} \le M_{0,10} \le N_{2,10} \le N_{1,10} \le_3 N_{0,10}.
	\end{align*}
\end{conj}

In this paper we prove this conjecture.

\begin{thm}\label{T:main}
	\Cref{CONJ:rankcrankconj} is true.
\end{thm}

The paper is organized as follows: In \Cref{S:Genfunc} we state generating functions for our combinatorial objects. In \Cref{S:Approxcranks} (resp. \Cref{S:Approxranks}) we determine the asymptotic behavior of the crank (resp. rank) with an explicit error term. \Cref{S:Proof} is devoted to the proof of \Cref{T:main}.

\section*{Acknowledgments}

The authors have received funding from the European Research Council (ERC) under the European Union’s Horizon
2020 research and innovation programme (grant agreement No. 101001179). The authors would like to thank Caner Nazaroglu for helpful discussions and great suggestions on Mathematica implementation.

\section{Generating functions}\label{S:Genfunc}

In this section, using orthogonality of the roots of unity, we write the generating function for the crank (resp. rank) in an arithmetic progression in terms of the crank (resp. rank) function specialized at roots of unity. We consider the generating functions
\begin{align*}
	C(\z;q) &:= \sum_{\substack{n\ge0\\m\in\Z}} M(m,n) \z^m q^n = \frac{(q;q)_\infty}{(\z q;q)_\infty\lrb{\z^{-1}q;q}_\infty},\\
	R(\z;q) &:= \sum_{\substack{n\ge0\\m\in\Z}} N(m,n) \z^m q^n = \sum_{n=1}^\infty \frac{q^{n^2}}{(\z q;q)_\infty\lrb{\z^{-1}q;q}_\infty},
\end{align*}
where $(a;q)_n:=\prod_{j=0}^{n-1}(1-aq^j)$, $a\in\C$, $n\in\N_0\cup\{\infty\}$. Using orthogonality of roots of units ($\z_c:=e^\frac{2\pi i}{c}$), we have
\[
	\sum_{n=0}^\infty M(a,c;n) q^n = \frac1c \sum_{n=0}^\infty p(n) q^n + \frac1c \sum_{j=1}^{c-1} \z_c^{-aj} C\left(\z_c^j;q\right),
\]
and
\[
	\sum_{n=0}^\infty N(a,c;n) q^n = \frac1c \sum_{n=0}^\infty p(n) q^n + \frac1c \sum_{j=1}^{c-1} \z_c^{-aj} R\left(\z_c^j;q\right).
\]
Now assume that $c$ is odd. Noting that
\[
	C\left(\z_c^{-j};q\right) = C\left(\z_c^j;q\right),\qquad R\left(\z_c^{-j};q\right) = R\left(\z_c^j;q\right),
\]
we may write
\begin{align*}
	\sum_{n=0}^\infty M(a,c;n)q^n &= \frac1c\sum_{n=0}^\infty p(n)q^n + \frac2c\sum_{j=1}^\frac{c-1}{2} \cos\pa{2\pi aj}{c}C\left(\z_c^j;q\right),\\
	\sum_{n=0}^\infty N(a,c;n)q^n &= \frac1c\sum_{n=0}^\infty p(n)q^n + \frac2c\sum_{j=1}^\frac{c-1}{2} \cos\pa{2\pi aj}{c}R\left(\z_c^j;q\right).
\end{align*}
To prove \Cref{T:main}, let
\begin{align*}
	\CC_a^{[1]}(n) &:= M(a,11;n) - \frac{p(n)}{11},\qquad \CC_a^{[2]}(n) := N(a,11;n) - \frac{p(n)}{11},\\
	\CC_{a_1,a_2}^{[3]}(n) &:= N(a_1,11;n) - N(a_2,11;n),\qquad \CC_{a_1,a_2}^{[4]}(n) := N(a_1,11;n) - M(a_2,11;n),\\
	\CC_{a_1,a_2}^{[5]}(n) &:= M(a_1,11;n) - M(a_2,11;n).
\end{align*}
Write
\[
	C\left(\z_{11}^j;q\right) =: \sum_{n\ge0} A_j(n)q^n,\qquad R\left(\z_{11}^j;q\right) =: \sum_{n\ge0} B_j(n)q^n.
\]
Then we have
\begin{align*}
	\CC_a^{[1]}(n) &= \frac{2}{11}\sum_{j=1}^5 \cos\pa{2\pi aj}{11}A_j(n),\qquad \CC_a^{[2]}(n) = \frac{2}{11}\sum_{j=1}^5 \cos\pa{2\pi aj}{11}B_j(n),\\
	\CC_{a_1,a_2}^{[3]}(n) &= \frac{2}{11}\sum_{j=1}^5 \left(\cos\pa{2\pi a_1j}{11}-\cos\pa{2\pi a_2j}{11}\right)B_j(n),\\
	\CC_{a_1,a_2}^{[4]}(n) &= \frac{2}{11}\sum_{j=1}^5 \left(\cos\pa{2\pi a_1j}{11}B_j(n)-\cos\pa{2\pi a_2j}{11}A_j(n)\right),\\
	\CC_{a_1,a_2}^{[5]}(n) &= \frac{2}{11}\sum_{j=1}^5 \left(\cos\pa{2\pi a_1j}{11}-\cos\pa{2\pi a_2j}{11}\right)A_j(n).
\end{align*}

\section{Approximating cranks}\label{S:Approxcranks}

The goal of this section is to approximate $M(j,11;n)$ with an explicit error. For $h,k$ coprime, let $h'$ be a solution to $hh'\equiv-1\Pmod k$ if $k$ is odd and to $hh'\equiv-1\Pmod{2k}$ if $k$ is even. For $c\mid k$, define the Kloosterman sum
\[
	B_{j,c,k}^*(n,m) := (-1)^{jk+1}\sin\pa{\pi j}{c}\sum_{h\Pmod k^*} \frac{\w_{h,k}}{\sin\pa{\pi jh'}{c}}e^{-\frac{\pi ij^2kh'}{c^2}}e^{\frac{2\pi i}{k}\left(nh+mh'\right)},
\]
where the sum runs over all $h$ modulo $k$ that are coprime to $k$ and we have (see \cite[equation (5.2.4)]{A})
\begin{equation*}
	\w_{h,k} :=
	\begin{cases}
		\pa{-k}{h}e^{-\pi i\left(\frac14(2-hk-h)+\frac{1}{12}\left(k-\frac1k\right)\left(2h-h'+h^2h'\right)\right)} & \text{if }h\text{ is odd},\\
		\pa{-h}{k}e^{-\pi i\left(\frac14(k-1)+\frac{1}{12}\left(k-\frac1k\right)\left(2h-h'+h^2h'\right)\right)} & \text{if }k\text{ is odd},
	\end{cases}
\end{equation*}
Here $(\frac\cdot\cdot)$ denotes the Kronecker symbol. Moreover, let $\d_S:=1$ if a statement $S$ is true and $\d_S:=0$ otherwise.

\begin{lem}\label{L:MCE}
	We have, for $1\le j\le5$,
	\[
		M(j,11;n) = M_j(n) + \sin\pa{\pi j}{11}\CE_j^{[1]}(n),
	\]
	where
	\[
		M_j(n) := \frac{4\sqrt3}{\sqrt{24n-1}} \sinh\pa{\pi\sqrt{24n-1}}{66} \left(\frac{iB_{j,11,11}^*(-n,0)}{\sqrt{11}} + 2\sin\pa{\pi}{11}\d_{j\equiv1\Pmod{11}}\right),
	\]
	and where $\CE_j^{[1]}(n)$ may be bound against 
	\begin{equation*}
		\left|\CE_j^{[1]}(n)\right|
		\le \CE_1(n) := \frac{184}{11\sqrt3}n^\frac34\sinh\pa{\pi\sqrt{24n-1}}{132} + 77491n^\frac14 + 1324.9n^\frac38 + 36620.2.
	\end{equation*}
\end{lem}

\subsection{The asymptotic formula}

From Theorem 1.1 of \cite{ZR} it is not hard to conclude that, for $\e>0$, we have
\begin{align}\label{E:exactcrank}
	&A_{j,11}(n) = \frac{4\sqrt3i}{\sqrt{24n-1}}\sum_{\substack{1\le k\le\sqrt n\\11\mid k}} \frac{B_{j,11,k}^*(-n,0)}{\sqrt k}\sinh\pa{\pi\sqrt{24n-1}}{6k}\\
	\nonumber
	&\hspace{.5cm}+ \frac{8\sqrt3\sin\pa{\pi j}{11}}{\sqrt{24n-1}}\sum_{\substack{1\le k\le\sqrt n\\11\nmid k\\jk\equiv\pm1\Pmod{11}}} \frac{D_{j,11,k}\left(-n,m_{j,11,k,0}^+\right)}{\sqrt k}\sinh\pa{\pi\sqrt{24n-1}}{66k} + O\left(n^\e\right),
\end{align}
with $m_{j,11,k,0}^+\in\Z$ and

\[
	D_{j,c,k}(n,m) := (-1)^{jk+\ell}\sum_{h\Pmod k^*} \w_{h,k}e^{\frac{2\pi i}{k}\left(nh+mh'\right)},
\]
where $1\le\ell\le{10}$ is the unique solution to $\ell\equiv kj\Pmod {11}$.

\subsection{The main term}

The term $k=11$ from the first sum and the term $k=1$ (if it occurs) from the second sum in \eqref{E:exactcrank} have the same size, and their sum is equal to
\begin{multline*}
	\frac{4\sqrt3i}{\sqrt{24n-1}}\frac{B_{j,11,11}^*(-n,0)}{\sqrt{11}}\sinh\pa{\pi\sqrt{24n-1}}{66}\\
	+ \frac{8\sqrt3\sin\pa{\pi j}{11}}{\sqrt{24n-1}}\d_{j\equiv1\Pmod{11}}D_{j,11,1}\left(-n,m_{j,11,1,0}^+\right)\sinh\pa{\pi\sqrt{24n-1}}{66}.
\end{multline*}
The main term in \Cref{L:MCE} follows by evaluating
\begin{equation*}
	D_{j,11,1}\left(-n,m_{j,11,1,0}^+\right) = (-1)^{j+\ell} = 1.
\end{equation*}

\subsection{Error bounds}

In this subsection we absorb the remaining terms into the error term. We first bound the contribution from the first term in \eqref{E:exactcrank} with $k\ge22$. This can be estimated against
\begin{align*}
	&\frac{4\sqrt3}{\sqrt{24n-1}}\sin\pa{\pi j}{11}\sum_{\substack{22\le k\le\sqrt n\\11\mid k}} \frac{1}{\sqrt k}\sum_{h\Pmod k^*} \frac{1}{\left|\sin\pa{\pi jh'}{11}\right|}\sinh\pa{\pi\sqrt{24n-1}}{6k}\\
	&\hspace{.5cm}\le \frac{4\sqrt3 \sin\left(\frac{\pi j}{11}\right) \sinh\pa{\pi\sqrt{24n-1}}{6\cdot22}}{\sqrt{24n-1}\sin\pa{\pi}{11}}\sum_{\substack{22\le k\le\sqrt n\\11\mid k}} \sqrt k \le \frac{4\sqrt3 \sin\left(\frac{\pi j}{11}\right) \sinh\pa{\pi\sqrt{24n-1}}{132}}{\sqrt{24n-1}\sin\pa{\pi}{11}}\frac{2}{33} n^\frac34\\
	&\hspace{.5cm}\le \frac{8\sin\pa{\pi j}{11}n^\frac34}{11\sqrt3\sqrt{24n-1}}\sinh\pa{\pi\sqrt{24n-1}}{132},
\end{align*}
where, in the penultimate step, we use that
\begin{equation}\label{E:sumint}
	\sum_{\substack{1\le k\le\sqrt n\\c\mid k}} k^\frac12 \le \frac{2}{3c}n^\frac34.
\end{equation}
The contribution from the second sum in \eqref{E:exactcrank} with $k\ge2$ can be bound against, using \eqref{E:sumint} with $c=1$,
\[
	\frac{8\sqrt3 \sin\left(\frac{\pi j}{11}\right)}{\sqrt{24n-1}}\sinh\pa{\pi\sqrt{24n-1}}{132}\sum_{1\le k\le\sqrt n} \sqrt k \le \frac{16\sin\left(\frac{\pi j}{11}\right)}{\sqrt3\sqrt{24n-1}}n^\frac34 \sinh\pa{\pi\sqrt{24n-1}}{132}.
\]

To make the error term in \eqref{E:exactcrank} explicit we follow the proof of \cite{ZR} (which follows the proof of \cite{Br}) using the notation in there. From \cite{ZR} the error that is obtained by bounding the non-principal part of $\Sigma_1$ is
\[
	|S_\err| \le 2e^{2\pi+\frac{\pi}{24}}\sin\tpa{\pi j}{11}\left(c_2+2\left(1+\cos\tpa{\pi}{11}\right)c_1(1+c_2)\right)\sum_{\substack{1\le k\le\sqrt n\\11\mid k}} k^{-\frac32}\sum_{\substack{1\le h\le k\\\gcd(h,k)=1}} \frac{1}{\left|\sin\pa{\pi h}{11}\right|},
\]
where
\[
	c_1 := \sum_{m=1}^\infty \frac{e^{-\frac{\pi m(m+1)}{2}}}{1-e^{-\pi m}},\qquad c_2 := \sum_{n=1}^\infty p(n)e^{-\pi n}.
\]
We bound
\[
	c_1 \le \frac{1}{1-e^{-\pi}}\sum_{m=1}^\infty e^{-\pi m} = \frac{e^{-\pi}}{\left(1-e^{-\pi}\right)^2}.
\]
Moreover, using that $p(n)<\exp(\pi\sqrt\frac{2n}{3})$, we estimate
\[
	p(n)e^{-\pi n} < e^{\pi\sqrt\frac{2n}{3}-\pi n} \le e^{-\frac{\pi n}{2}} \qquad\text{for }n\ge3.
\]
Thus we get
\begin{equation}\label{E:Pbound}
	c_2 \le p(1)e^{-\pi} + p(2)e^{-2\pi} + \sum_{n=3}^\infty e^{-\frac{\pi n}{2}} = e^{-\pi} + 2e^{-2\pi} + \frac{e^{-\frac{3\pi}{2}}}{1-e^{-\frac\pi2}}.
\end{equation}
Moreover, we have
\[
	\sum_{\substack{1\le k\le\sqrt n\\11\mid k}} k^{-\frac32}\sum_{\substack{1\le h\le k\\\gcd(h,k)=1}} \frac{1}{\left|\sin\pa{\pi h}{11}\right|} \le \frac{1}{\sin\pa{\pi}{11}}\sum_{\substack{1\le k\le\sqrt n\\11\mid k}} k^{-\frac12} \le \frac{2n^\frac14}{11\sin\pa{\pi}{11}},
\]
using that
\begin{equation}\label{E:int}
	\sum_{\substack{1\le k\le\sqrt n\\11\mid k}} k^{-\frac12} = 11^{-\frac12}\sum_{1\le k\le\frac{\sqrt n}{11}} \frac{1}{k^\frac12} \le 11^{-\frac12}\int_0^\frac{\sqrt n}{11} x^{-\frac12} dx = 11^{-\frac12}2\sqrt\frac{\sqrt n}{11} = \frac{2n^\frac14}{11}.
\end{equation}
Thus, we get
\[
	|S_\err| \le 200.2 \sin\pa{\pi j}{11}n^\frac14.
\]

We next estimate the contribution from the non-principal part in $\Sigma_2$. From \cite{ZR} this may be bound against
\begin{equation}\label{E:Terr}
	|T_\err| \le 16e^{2\pi}f(11)\sin\pa{\pi j}{11}n^\frac14,
\end{equation}
where
\[
	f(c) := \frac{1+c_2e^{\pi\d_0{ (c)}}}{1-e^{-\frac{\pi}{c}}} + e^{\pi\d_0{ (c)}}c_1(1+c_2) + \frac12e^{\pi\d_0{ (c)}}(c_2+1)c_3,
\]
where
\begin{equation*}
	\d_0(c) := \frac{1}{2 c^2} + \frac{1}{24} - \frac{1}{2 c}
\end{equation*}
and thus $\d_0(11) = \frac{1}{2904}$.
Moreover we have
\begin{equation*}
	c_3 := \sum_{m=2}^\infty \frac{e^{-\frac{\pi m(m+1)}{2}}}{1-e^{\pi-\pi m}} \le \frac{1}{1-e^{-\pi}}\sum_{m=2}^\infty e^{-\frac{3\pi m}{2}} = \frac{e^{-3\pi}}{\left(1-e^{-\pi}\right)\left(1-e^{-\frac{3\pi}{2}}\right)}.
\end{equation*}
Plugging these bounds into \eqref{E:Terr} we obtain 
\[
	|T_\err| \le 36928.5\sin\pa{\pi j}{11}n^\frac14.
\]

Next the error obtained by symmetrizing the first main contribution may be bound again, using \eqref{E:int},
\begin{align*}
	\left|S_\err^{[1]}\right| &\le \sin\pa{\pi j}{11}e^{2\pi+\frac{\pi}{12}}n^\frac12\frac{2}{\Flo{\sqrt n} + 1}\sum_{\substack{0\le h<k\le\sqrt n\\\gcd(h,k)=1\\11\mid k}} \frac{1}{k^\frac32}\frac{1}{\left|\sin\pa{\pi h}{11}\right|}\\
	&\le 2\sin\pa{\pi j}{11}e^{2\pi+\frac{\pi}{12}}\frac{1}{\sin\pa{\pi}{11}}\frac{2}{11}n^\frac14 \le 898.1\sin\pa{\pi j}{11}n^\frac14.
\end{align*}

The error obtained by symmetrizing the second main term contribution may be estimated against
\[
	\left|S_\err^{[2]}\right| \le 32e^{2\pi}\sin\pa{\pi j}{11}\frac{e^{2\pi\d_0}}{1-e^{-\frac{2\pi}{11}}}n^\frac14 \le 39464.1\sin\pa{\pi j}{11}n^\frac14.
\]
Next we have the errors obtained by integrating against the smaller arc
\begin{align*}
	\frac{4\left(\frac43+2^\frac54\right)\sin\pa{\pi j}{11}(1+\log(5))e^{2\pi+\frac{\pi}{12}}}{11\pi\left(1-\frac{\pi^2}{24}\right)}n^\frac38 &\le 1324.9\sin\pa{\pi j}{11}n^\frac38,\\
	8\left(\frac43+2^\frac54\right)\sin\pa{\pi j}{11}\frac{e^{2\pi\d_0+2\pi}}{1-e^{-\frac{2\pi}{11}}} &\le 36620.2\sin\pa{\pi j}{11}.
\end{align*}
Combining gives that overall error can be bound against
\begin{multline*}
	\sin\pa{\pi j}{11}\left(\frac{8n^\frac34}{11\sqrt3\sqrt{24n-1}}\sinh\pa{\pi\sqrt{24n-1}}{132} + \frac{16n^\frac34}{\sqrt3\sqrt{24n-1}}\sinh\pa{\pi\sqrt{24n-1}}{132}\right.\\
	\left.+ 200.2n^\frac14 + 36928.5n^\frac14 + 898.1n^\frac14 + 39464.1n^\frac14 + 1324.9n^\frac38 + 36620.2{\vphantom{\frac{n^\frac34}{\sqrt3}}}\right).
\end{multline*}
Simplifying gives the claim.

\section{Approximating ranks}\label{S:Approxranks}

The goal of this section is to show the following approximation of $B_j(n)$.

\begin{lem}\label{L:rank}
	We have, for $1\le j\le5$,
	\[
		B_j(n) = \frac{8\sqrt3\sin\pa{\pi}{11}}{\sqrt{24n-1}}\d_{j\equiv1\Pmod{11}} \sinh\pa{5\pi\sqrt{24n-1}}{66} + \sin\pa{\pi j}{11}\CE_j^{[2]}(n),
	\]
	where $\CE_j^{[2]}$ may be bound against 
	\[
		\left|\CE_j^{[2]}(n)\right| \le \CE_2(n) := \frac{10.74n^\frac34}{\sqrt{24n-1}}\sinh\pa{\pi\sqrt{24n-1}}{132} + 18092.8n^\frac14 + 40707.2n^\frac54.
	\]
\end{lem}

\subsection{The asymptotic formula}

Using Theorem 1.1 of \cite{Br} it is not hard to show that
\begin{align}\label{E:Bn}
	&B_j(n) = \frac{4\sqrt3i}{\sqrt{24n-1}}\sum_{\substack{1\le k\le\sqrt n\\11\mid k}} \frac{B_{j,11,k}(-n,0)}{\sqrt k}\sinh\pa{\pi\sqrt{24n-1}}{6k}\\
	\nonumber
	&\hspace{1.6cm}+ \frac{8\sqrt3\sin\pa{\pi j}{11}}{\sqrt{24n-1}}\sum_{\substack{1\le k\le\sqrt n\\11\nmid k\\jk\equiv\pm1\Pmod{11}}} \frac{D_{j,11,k}(-n,m_{j,11,k,0})}{\sqrt k}\sinh\pa{5\pi\sqrt{24n-1}}{66k} + O\left(n^\e\right).
\end{align}
Here $m_{j,11,k,0}\in\Z$, and for $c\mid k$, we have 
\begin{align*}
	B_{j,c,k}(n,m) &:= (-1)^{jk+1}\sin\pa{\pi j}{c}\sum_{h\Pmod k^*} \frac{\w_{h,k}}{\sin\pa{\pi jh'}{c}}e^{-\frac{3\pi ij^2kh'}{c}}e^{\frac{2\pi i}{k}\left(nh+mh'\right)},\\
	\intertext{and}
	D_{j,c,k}(n,m) &:= (-1)^{jk+\ell}\sum_{h\Pmod k^*} \w_{h,k}e^{\frac{2\pi i}{k}\left(nh+mh'\right)},
\end{align*}
where, just like in the case of cranks, $1\le\ell\le10$ is the unique solution to $jk\equiv\ell\Pmod{11}$.

\subsection{The dominant term}

The main contribution in \Cref{L:rank} comes from the term $k=1$, which only occurs if $\ell=1$ and thus $j=1$.

\subsection{Error bounds}

As in the case of the crank, we first bound the first sum in \eqref{E:Bn} against
\begin{equation*}
	\frac{4\sqrt3}{\sqrt{24n-1}} \sum_{\substack{1\le k\le\sqrt n\\11\mid k}} \hspace{-.1cm} \frac{|B_{j,11,k}(-n,0)|}{\sqrt k}\sinh\hspace{-.05cm}\pa{\pi\sqrt{24n-1}}{6k}
	\le \frac{1.5\sin\left(\frac{\pi j}{11}\right)n^\frac34}{\sqrt{24n-1}}\sinh\hspace{-.05cm}\pa{\pi\sqrt{24n-1}}{66},
\end{equation*}
using again \eqref{E:sumint}.

Similarly, the contribution from $k\ne1$ from the second sum can be bound against
\[
	\frac{16\sin\pa{\pi j}{11}}{\sqrt3\sqrt{24n-1}}n^\frac34\sinh\pa{5\pi\sqrt{24n-1}}{132}.
\]

We now look at the $O$-terms and make these explicit. We follow the proof and the notation of \cite{Br}. We bound the Kloosterman sums trivially. The contribution from $S_2$ may be estimated against
\begin{align*}
	|S_2| &\le \sin\pa{\pi j}{11}\sum_{\substack{0\le h<k\le\sqrt n\\11\mid k\\\gcd(h,k)=1}} \frac{1}{\left|\sin\pa{\pi jh'}{11}\right|}\frac{2}{k\left(\Flo{\sqrt n}+1\right)}\frac{\sqrt n}{\sqrt k}e^{\frac{2\pi}{k}\frac kn\left(n-\frac{1}{24}\right)} \sum_{r\ge1} |a(r)|e^{-\frac{2\pi}{k}\left(r-\frac{1}{24}\right)\frac k2}\\
	&\le \frac{2e^{2\pi}\sqrt n\sin\left(\frac{\pi j}{11}\right)}{\left(\Flo{\sqrt n}+1\right)\sin\pa{\pi}{11}}\sum_{r\ge1} |a(r)|e^{-\pi\left(r-\frac{1}{24}\right)}\sum_{\substack{1\le k\le\sqrt n\\11\mid k}} \frac{1}{k^\frac12}.
\end{align*}
To bound the sum on $k$, we use \eqref{E:int}. For the sum on $r$, we require the estimate from Section 5 of \cite{Br}
\[
	\sum_{r\ge1} |a(r)|e^{-\pi r} \le 0.06.
\]
Combining gives that
\[
	|S_2| \le \frac{2\cdot0.06\sin\pa{\pi j}{11}\cdot2}{11\sin\pa{\pi}{11}}e^{2\pi+\frac{\pi}{24}}n^\frac14 \le 47.3\sin\pa{\pi j}{11}n^\frac14.
\]

We next bound the non-principal part from $\sum_2$ against 
\begin{equation}\label{E:nonp}
	4\sin\pa{\pi j}{11}\sum_{\substack{0\le h<k\le\sqrt n\\\gcd(h,k)=1\\11\nmid k}} \frac{2}{k\left(\Flo{\sqrt n}+1\right)}\frac{\sqrt n}{\sqrt k}e^{\frac{2\pi}{n}\left(n-\frac{1}{24}\right)}\sum_{r\ge1} |b(r)|e^{-\frac{\pi r}{12k\cdot11^2}\frac k2},
\end{equation}
where $b(r)$ is defined via
\begin{equation}\label{E:b-r-def}
	e^{-\frac{2\pi i h' sa}{c} - \frac{3 \pi i h' a^{2}k_{1}}{cc_{1}} + \frac{6 \pi i h' \ell a}{cc_{1}} + \frac{\pi}{12 k z}} q_{1}^{\frac{s(\ell,11)\ell}{c_{1}} - \frac{3 \ell^{2}}{2c_{1}^{2}}} N\left(ah', \frac{\ell c}{c_{1}}, c; q_{1}\right)
	=:\sum_{r \ge r_{0}} b(r)e^{\frac{2\pi i m_{r}h'}{k}}e^{-\frac{\pi i r}{12k c^{2}z}},
\end{equation}
where $s(\ell,c)$ is defined in (1.9) of \cite{Br}.
We first show that
\begin{equation}\label{E:boundb}
\sum_{r \ge 1} \lvert b(r) \rvert e^{-\tfrac{\pi r}{2904}} \leq 2.8,
\end{equation}
By taking $k=1, c=1, h'=0$ in \eqref{E:b-r-def} we have, for $q \in \R$,
\begin{equation*}
q^{-\frac{1}{24} + \frac{s(\ell,11)\ell}{11} - \frac{3\ell^{2}}{2 \cdot 11^{2}}} N\left(0,\ell,11;q\right) = \sum_{r \ge 0} b(r) q^{\frac{r}{24 \cdot 11^{2}}}.
\end{equation*}
By definition
\begin{multline*}
	N(0,\ell,11;q)\\
	= \frac{i}{2\left(q;q\right)_{\infty}} \left(\sum_{m \ge 0} \frac{(-1)^{m}q^{\frac{m}{2}(3m + 1)+ms(\ell,11) + \frac{\ell}{22}}}{1 - q^{m + \frac{\ell}{11}}}
	-\sum_{m \ge 1} \frac{(-1)^{m}q^{\frac{m}{2}(3m + 1) - ms(\ell,11) - \frac{\ell}{22}}}{1 - q^{m - \frac{\ell}{11}}}\right).
\end{multline*}
Thus, we have
\begin{multline*}
	\sum_{r\ge0} |b(r)|q^{\frac{r}{24\cdot11^2}+\frac{1}{24}-\frac{s(\ell,11)\ell}{11} + \frac{3\ell^2}{2\cdot11^2}}\\
	= \frac{q^\frac{\ell}{22}}{2}P(q)\sum_{m\ge0} \frac{q^{\frac{m(3m+1)}{2}+ms(\ell,11)}}{1-q^{m+\frac{\ell}{11}}} + \frac{q^{-\frac{\ell}{22}}}{2}P(q)\sum_{m\ge1} \frac{q^{\frac{m(3m+1)}{2}-ms(\ell,11)}}{1-q^{m-\frac{\ell}{11}}},
\end{multline*}
where $P(q):=\frac{1}{\lrb{q;q}_{\infty}}=\sum_{n\ge0} p(n)q^n$. We bound the first sum, for $0<q<1$,
\begin{equation*}
	\sum_{m \ge 0} \frac{q^{\frac{m(3m+1)}2+ms(\ell,11)}}{1-q^{m+\frac\ell{11}}} \leq \frac1{\left(1-q^{\frac\ell{11}}\right)\left(1-q^{\frac12+s(\ell,11)}\right)}.
\end{equation*}
For the second sum, we estimate
\begin{equation*}
	\sum_{m \ge 1} \frac{q^{\frac{m(3m+1)}2-ms(\ell,11)}}{1-q^{m-\frac\ell{11}}} \leq \frac{q^{2-s(\ell,11)}}{1-q^{1-\frac\ell{11}}} + \frac{q^{2\left(\frac72-s(\ell,11)\right)}}{\left(1-q^{2-\frac\ell{11}}\right)\left(1-q^{\frac72-s(\ell,11)}\right)}.
\end{equation*}
Thus, for $q=e^{-\pi}$ and using the bound \eqref{E:Pbound}, we get
\begin{align*}
	&\sum_{r\ge0} |b(r)|e^{-\frac{\pi r}{2904}} \le e^{-\pi\left(-\frac{1}{24}+\frac{s(\ell,11)\ell}{11}-\frac{3\ell^2}{2\cdot11^2}\right)} \frac12\left(1+e^{-\pi}+2e^{-2\pi}+\frac{e^{-\frac{3\pi}{2}}}{1-e^{-\frac\pi2}}\right)\\
	&\hspace{1.5cm}\times \left(\frac{e^{-\frac{\pi\ell}{22}}}{\left(1-e^{-\frac{\pi\ell}{11}}\right)\left(1-e^{-\pi\left(\frac12+s(\ell,11)\right)}\right)}\right.\\
	&\hspace{3cm}+ \left.e^{\frac{\pi\ell}{22}}\left(\frac{e^{-\pi(2-s(\ell,11))}}{1-e^{-\pi\left(1-\frac\ell{11}\right)}}+\frac{e^{-2\pi\left(\frac72-s(\ell,11)\right)}}{\left(1-e^{-\pi\left(2-\frac\ell{11}\right)}\right)\left(1-e^{-\pi\left(\frac72-s(\ell,11)\right)}\right)}\right)\right).
\end{align*}
We now see that this contribution is maximized for $\ell=1$ (in which case $s(\ell,11)=0$) and in this case can be bound against $\leq2.8$. This gives \eqref{E:boundb}. 

Thus \eqref{E:nonp} may be estimated by 
\begin{equation*}
	\frac{4\sin\pa{\pi j}{11}2\sqrt n}{\Flo{\sqrt n}+1}e^{2\pi}\cdot2.8\sum_{1\le k\le\sqrt n} \hspace{-0.1cm} \frac{1}{\sqrt k} \le 4\sin\left(\frac{\pi j}{11}\right)2e^{2\pi}\cdot2.8\cdot2n^\frac14 \le 8567.9\sin\left(\frac{\pi j}{11}\right)n^\frac14,
\end{equation*}
bounding the sum on $k$ against
\begin{equation}\label{E:int2}
	\int_0^{\sqrt n} x^{-\frac12} dx = 2n^\frac14.
\end{equation}

Next the contribution from $S_{12}$ gives ($S_{13}$ yields a contribution of exactly the same size)
\begin{align*}
	|S_{12}| &\le \sin\pa{\pi j}{11}\sum_{\substack{0\le h<k\le\sqrt n\\11\mid k\\\gcd(h,k)=1}} \frac{1}{\left|\sin\pa{\pi jh'}{11}\right|}\frac{2}{k\left(\Flo{\sqrt n}+1\right)}\frac{\sqrt n}{\sqrt k}e^{\frac{2\pi}{k}\frac kn\left(n-\frac{1}{24}\right)}e^{\frac{\pi}{12k}k}\\
	&\le \frac{2\sin\left(\frac{\pi j}{11}\right)e^{2\pi+\frac{\pi}{12}}}{\sin\pa{\pi}{11}}\sum_{\substack{1\le k\le\sqrt n\\11\mid k}} \frac{1}{\sqrt k} \le \frac{2\sin\left(\frac{\pi j}{11}\right)e^{2\pi+\frac{\pi}{12}}}{\sin\pa{\pi}{11}}\frac{2}{11}n^\frac14 \le 898.1\sin\pa{\pi j}{11}n^\frac14,
\end{align*}
where we use \eqref{E:int} to estimate the sum on $k$. Next the corresponding contribution from $\sum_2$ may be bound against
\[
	\sin\pa{\pi j}{11}\sum_{\substack{0\le h<k\le\sqrt n\\\gcd(h,k)=1\\11\nmid k\\r\ge0\\\d_{11,k,r}>0}} \frac{2}{k\left(\Flo{\sqrt n}+1\right)}\frac{\sqrt n}{\sqrt k}e^{\frac{2\pi}{n}\left(n-\frac{1}{24}\right)}e^{\frac{2\pi\d_{11,k,r}}{k}k},
\]
where $\d_{c,k,r}$ is defined in (1.13) of \cite{Br}. We determine for which values $\d_{11,k,r}>0$, namely $jk\equiv\pm1\Pmod{11}$, $r=0$, and $\d_{11,k,0}=\frac{25}{2904}$ in both cases. Thus we bound the above against, using \eqref{E:int2},
\begin{align*}
	2\sin\left(\frac{\pi j}{11}\right)e^{2\pi\frac{25}{2904}+2\pi}\sum_{\substack{1\le k\le\sqrt n\\11\nmid k\\jk\equiv\pm1\Pmod{11}}} \frac{1}{\sqrt k} \le 4\sin\left(\frac{\pi j}{11}\right)n^\frac14e^{2\pi\frac{25}{2904}+2\pi}
	\le 2261.1\sin\left(\frac{\pi j}{11}\right)n^\frac14.
\end{align*}

On the smaller arcs we have the same estimates as on $S_{12}$ and we obtain a contribution of the sine.

Finally $\Sigma_3$ may be bounded against
\[
	\left|\ssum_3\right| \le 2\sin\pa{\pi j}{11}^2\sum_{\substack{0\le h<k\le\sqrt n\\\gcd(h,k)=1}} \frac1k\sum_{\nu=0}^{k-1} \left|\int_{-\vth_{h,k}'}^{\vth_{h,k}''} e^{\frac{2\pi z}{k}\left(n-\frac{1}{24}\right)}z^\frac12I_{j,11,k,\nu}(z) d\Phi\right|,
\]
where $z:=\frac kn+k\Phi i$ and where the Mordell integral $I_{j,11,k,\nu}(z)$ is defined in \cite{Br}. We next bound the Mordell integrals, making Lemma 3.1 of \cite{Br} explicit as ($H_{j,c}$ is defined in \cite{Br})
\begin{multline*}
	\left|z^\frac12 I_{j,11,k,\nu}(z)\right| \le \frac k{\pi|z|^\frac12} \int_\R \left|e^{-\frac{3kt^2}{\pi z}}\right| \left|H_{j,c}\left(\frac{\pi i\nu}k - \frac{\pi i}{6k} - t\right)\right| dt\\
	= \frac k{\pi|z|^\frac12} \int_\R e^{-\frac{3k}\pi \re\left(\frac1z\right) t^2} \left|\frac{\cosh\left(\frac{\pi i\nu}k - \frac{\pi i}{6k} - t\right)}{\sinh\left(\frac{\pi i\nu}k - \frac{\pi i}{6k} - t + \frac{\pi ij}{11}\right) \sinh\left(\frac{\pi i\nu}k - \frac{\pi i}{6k} - t - \frac{\pi ij}{11}\right)}\right| dt.
\end{multline*}
Now, for $a,b\in\R$,
\[
	\left|\cosh\left(\frac{\pi i\nu}{k}-\frac{\pi i}{6k}-t\right)\right| \le \cosh(t),\qquad \left|\sinh(a+ib)\right|^2 \ge
	\begin{cases}
		\frac{e^{2|a|}}{6} & \text{if } |a|\ge1,\\
		\sin(b)^2 & \text{if } |a|\le1.
	\end{cases}
\]
Thus, using that $\re(\frac1z)>\frac k2$
\begin{align*}
	\left|z^\frac12 I_{j,11,k,\nu}(z)\right| &\le \frac{k}{\pi|z|^\frac12}\left(6\int_{|t|\ge1} e^{-\frac{3t^2}{2\pi}}\cosh(t)e^{-2|t|} dt{\vphantom{\frac{e^{-\frac{t^2}{2}}}{\left|\sin\left(\frac{\pi j}{11}\right)\right|}}}\right.\\
	&\hspace{3.5cm}\left.+ \int_{|t|\leq1} \frac{e^{-\frac{3t^2}{2\pi}} \cosh(t)}{\left|\sin\left(\frac{\pi \nu}k - \frac\pi{6k} + \frac{\pi j}{11}\right)\right|\left|\sin\left(\frac{\pi \nu}k - \frac\pi{6k} - \frac{\pi j}{11}\right)\right|} dt\right)\\
	&\le \frac k{|z|^\frac12} \left(0.21 + \frac{0.64}{\left|\sin\left(\frac{\pi \nu}k - \frac\pi{6k} + \frac{\pi j}{11}\right)\right|\left|\sin\left(\frac{\pi \nu}k - \frac\pi{6k} - \frac{\pi j}{11}\right)\right|}\right).
\end{align*}
Now note that $\frac\nu k-\frac{1}{6k}\pm\frac{j}{11}\notin\Z$. Thus we bound
\begin{align*}
	\left|\sin\left(\frac{\pi\nu}{k}-\frac{\pi}{6k}+\frac{\pi j}{11}\right)\right|\left|\sin\left(\frac{\pi\nu}{k}-\frac{\pi}{6k}-\frac{\pi j}{11}\right)\right| &\ge \sin\pa{\pi}{66k}\sin\left(\frac{\pi}{66k}+\frac{\pi}{11}\right)\\
	&\ge \frac{\pi}{66k}\sin\pa{\pi}{11},
\end{align*}
using that for $0\le x\le\frac\pi2$, $\sin(x)\sin(x+\frac{\pi}{11})\ge x\sin(\frac{\pi}{11})$. This gives that
\[
	\left|z^\frac12 I_{j,11,k,\nu}(z)\right| \leq \frac k{|z|^\frac12} (0.21+47.73k) \leq \frac{48k^2}{|z|^\frac12}.
\]
Using $|z|\geq\frac kn$ we obtain 
\begin{align*}
	\left|\ssum_3\right| &\le 2\sin\pa{5\pi}{11}\sin\pa{\pi j}{11}48\sum_{\substack{0\le h<k\le\sqrt n\\\gcd(h,k)=1}} \frac1kk^2\pa kn^{-\frac12}\sum_{\nu=0}^{k-1} \frac{2}{k\left(\Flo{\sqrt n}+1\right)}e^{\frac{2\pi}{k}\frac kn\left(n-\frac{1}{24}\right)}\\
	&\le \frac{4\cdot48\sin\pa{5\pi}{11}\sin\pa{\pi j}{11}\sqrt ne^{2\pi}}{\Flo{\sqrt n}+1}\sum_{1\le k\le\sqrt n} k^\frac32 \le 40707.2\sin\pa{\pi j}{11}n^\frac54.
\end{align*}
Combining gives as error 
\begin{multline*}
	\sin\pa{\pi j}{11}\left(\frac{1.5n^\frac34}{\sqrt{24n-1}}\sinh\pa{\pi\sqrt{24n-1}}{66} + \frac{16}{\sqrt3\sqrt{24n-1}}n^\frac34\sinh\pa{5\pi\sqrt{24n-1}}{132}\right.\\
	\left.+ 47.3n^\frac14 + 8567.9n^\frac14 + 3\cdot898.1n^\frac14 + 3\cdot2261.1n^\frac14 + 40707.2n^\frac54{\vphantom{\frac{n^\frac34}{\sqrt{24n-1}}}}\right).
\end{multline*}
Simplifying gives the claim.

\section{Proof of \Cref{T:main}}\label{S:Proof}

We are now ready to prove \Cref{T:main}.

\subsection{$\mathcal C_a^{[1]}(n)$}

We need to show that
\[
	M_{a,d}(n) - \frac{p_{11}(n)}{11}
	\begin{cases}
		>0 & \text{if }(a,d)\in\CS_1,\\
		<0 & \text{if }(a,d)\in\CS_2,
	\end{cases}
\]
where
\begin{align*}
	\CS_1 &:= \{(0,0),(1,1),(2,2),(0,3),(0,4),(0,5),(1,7),(1,8),(1,9),(0,10)\},\\
	\CS_2 &:= \{(1,0),(2,1),(0,2),(1,3),(1,4),(2,5),(0,7),(0,8),(0,9),(3,10)\}.
\end{align*}
For $(j,d)\in\{(0,8),(1,8)\}$, we need to additionally assume that $n\geq41$. To show the claim, we require
\[
2\sum_{j=1}^5 \cos\pa{2\pi aj}{11}A_j(n)
\begin{cases}
	> 0 & \text{if }(a,d)\in\CS_1,\\
	< 0 & \text{if }(a,d)\in\CS_2,
\end{cases}
\]
where $n\equiv d\Pmod{11}$. A computer check gives that the claim holds for $n\le 1779$ (up to the above mentioned exceptions).

For $n\ge1780$, we use \Cref{L:MCE} so that the left-hand side equals
\begin{multline*}
	\tfrac{4\sqrt3}{\sqrt{24n-1}}\sinh\tpa{\pi\sqrt{24n-1}}{66} \left(\tfrac{2i}{\sqrt{11}}\sum_{j=1}^5 \cos\tpa{2\pi aj}{11}B_{j,11,11}^*(-n,0)
	+ 4\sin\tpa{\pi}{11}\cos\tpa{2\pi a}{11}\right)\\ +2\sum_{j=1}^5 \sin\tpa{\pi j}{11}\cos\tpa{2\pi aj}{11}\CE_j^{[1]}(n).
\end{multline*}
A numerical check gives that for $(a,d)\in\CS_1$, the first summand is positive and for $(a,d)\in\CS_2$, it is negative.

Let for $(a,d)\in\CS_1\cup\CS_2$
\[
	g_{a,d}(n) := 4\sqrt3\left(\frac{2i}{\sqrt{11}}\sum_{j=1}^5 \cos\pa{2\pi aj}{11}B_{j,11,11}^*(-11n-d,0)+4\sin\pa{\pi}{11}\cos\pa{2\pi a}{11}\right).
\]
We compute
\begin{align*}
		&iB_{j,11,11}^*(-m,0)\\
		&\hspace{.4cm}= 2(-1)^j\sin\pa{\pi j}{11}\left(-\frac{\cos\left(\frac{2\pi}{11}\left(m-\frac{j^2}{2}+10\right)\right)}{\sin\pa{\pi j}{11}}-\frac{\cos\left(\frac{2\pi}{11}\left(2m+5\frac{j^2}{2}+7\right)\right)}{\sin\pa{5\pi j}{11}}\right.\\
		&\hspace{.8cm}\left.- \frac{\cos\left(\frac{2\pi}{11}\left(3m-2j^2+2\right)\right)}{\sin\pa{4\pi j}{11}}-\frac{\cos\left(\frac{2\pi}{11}\left(4m-3\frac{j^2}{2}+2\right)\right)}{\sin\pa{3\pi j}{11}}+\frac{\cos\left(\frac{2\pi}{11}\left(5m+j^2+4\right)\right)}{\sin\pa{2\pi j}{11}}\right).
	\end{align*} 
A computer check gives that for $n\ge1780$, we have
\[
	\frac{2\sqrt{24(11n+d)-1}\CE_1(11n+d)}{\sinh\pa{\pi\sqrt{24(11n+d)-1}}{66}}\max_{(a,d)\in\CS_1\cup\CS_2} \frac{\sum_{j=1}^5 \left|\sin\pa{\pi j}{11}\cos\pa{2\pi aj}{11}\right|}{g_{a,d}(n)} < 1.
\]
This gives the claim.

\subsection{$\mathcal C_a^{[2]}(n)$}

We need to show that
\begin{multline*}
	N_{3,6}(n) \le \frac{p_{11}(n)}{11} \le N_{2,6}(n)\\
	\Leftrightarrow 2\sum_{j=0}^5 \cos\pa{6\pi j}{11}B_j( 11n+6) < 0,\quad 2\sum_{j=0}^5 \cos\pa{4\pi j}{11}B_j( 11n+6) > 0 \quad (n\equiv6\Pmod{11}).
\end{multline*}
A computer check gives that this is true for $n\le45$.

For $n\ge46$, we use \Cref{L:rank}. This gives for the sums of interest ($a\in\{2,3\}$)
\begin{align*}
	2\cos\pa{2\pi a}{11}\frac{8\sqrt3\sin\pa{\pi}{11}}{\sqrt{24 (11n+6)}-1}\sinh\pa{5\pi\sqrt{24 (11n+6)-1}}{66}\\
	 + 2\sum_{j=1}^5 \cos\pa{2\pi aj}{11}\sin\pa{\pi j}{11}\CE_j^{[2]}( 11n+6).
\end{align*}
So the sign of the main term is dictated by $\cos(\frac{2\pi a}{11})$ which is positive for $a=2$ and negative for $a=3$ as desired. 

To finish the claim, one can show that for $n\ge46$,
\begin{equation*}
	\frac{2\CE_2{ (11n+6)}\sum_{j=1}^5 \left|\cos\pa{2\pi aj}{11}\right|\sin\left(\frac{\pi j}{11}\right)}{2\left|\cos\pa{2\pi a}{11}\right|\frac{8\sqrt3\sin\pa{\pi}{11}}{\sqrt{24{ (11n+6)}-1}} \sinh\pa{5\pi\sqrt{24{ (11n+6)}-1}}{66}} < 1.
\end{equation*}

\subsection{$\mathcal C_{a_1,a_2}^{[3]}(n)$}

We need to show that
\[
	N_{5,d}(n) \le N_{4,d}(n) \le N_{3,d}(n),\qquad N_{2,d}(n) \le N_{1,d}(n) \le N_{0,d}(n).
\]
If $(d,a_1,a_2)\in\{(3,1,2),(5,3,4),(6,4,5),(6,0,1),(7,0,1)\}$, then we need to additionally assume that $n>11$, if $(d,a_1,a_2)\in\{(0,3,4),(3,0,1),(9,0,1)\}$, then we require $n>22$, and if $(d,a_1,a_2)\in\{(0,0,1),(2,0,1),(4,3,4),(4,0,1),(8,0,1),(10,0,1)\}$, then we let $n>33$. The claim follows if we prove that
\[
	N_{a,d}(n) - N_{a+1,d}(n) \ge 0 \text{ for } 0 \le a \le 4.
\]
Thus we need
\[
	\CC_{a,a+1}^{[3]}{ (11n+d)} \ge 0 \Leftrightarrow 2\sum_{j=1}^5 \left(\cos\pa{2\pi aj}{11}-\cos\pa{2\pi(a+1)j}{11}\right)B_j{ (11n+d)} \ge 0.
\]
A computer check gives that this is true for $n\le{46}$ (up to the above mentioned exception).

For $n\geq{ 47}$, we again use \Cref{L:rank}. The overall difference is
\begin{multline*}
	\tfrac{8\sqrt3}{\sqrt{24{(11n+d)}-1}}\sinh\left(\tfrac{5\pi}{66}\sqrt{24{ (11n+d)}-1}\right) 2\left(\cos\tpa{2\pi a}{11}-\cos\tpa{2\pi(a+1)}{11}\right)\sin\tpa{\pi}{11}\\
	+ 2\sum_{j=1}^5 \left|\cos\tpa{2\pi aj}{11}-\cos\tpa{2\pi(a+1)j}{11}\right|\sin\tpa{\pi j}{11}\CE_j^{[2]}{ (11n+d)}.
\end{multline*}
Now
\[
	\cos\pa{2\pi a}{11} - \cos\pa{2\pi(a+1)}{11} > 0.
\]

To finish the claim, we show that, for $0\le a\le4$, and $n\ge 47$,\begin{align*}
	\frac{2\CE_2{ (11n+d)}\sum_{j=1}^5 \left|\cos\pa{2\pi aj}{11}-\cos\pa{2\pi(a+1)j}{11}\right|\sin\pa{\pi j}{11}}{\frac{8\sqrt3}{\sqrt{24{ (11n+d)}-1}}\sinh\left(\frac{5\pi}{66}\sqrt{24{ (11n+d)}-1}\right) \cdot2 \sin\pa{\pi}{11}\left(\cos\pa{2\pi a}{11}-\cos\pa{2\pi(a+1)}{11}\right)} < 1.
\end{align*}

\subsection{$\mathcal C_{a_1,a_2}^{[4]}(n)$}

We need to show that
\[
	M_{a_1,d}(n) - N_{a_2,d}(n)
	\begin{cases}
		> 0 & \text{if } (a_1,a_2,d)\in \mathbb S_1,\\
		< 0 & \text{if } (a_1,a_2,d)\in \mathbb S_2,
	\end{cases}
\]
where
\begin{align*}
	\S_1 := \{&(1,3,0),(0,3,1),(0,3,2),(1,3,3),(1,3,4),(2,3,5),(0,3,7),(0,3,8),\\
	&(0,3,9),(3,3,10)\},\\
	\S_2 := \{&(0,2,0),(1,2,1),(2,2,2),(0,2,3),(0,2,4),(0,2,5),(1,2,7),(1,2,8),\\
	&(1,2,9),(0,2,10)\}.
\end{align*}
If $(a_1,a_2,d)\in\{(0,2,0),(0,3,1),(1,2,1),(2,2,2),(1,3,4),(0,2,4),(0,3,9),(3,3,10)\}$, then we need to additionally assume that $n\ge12$. A computer check shows that the claim holds for $n\le44$ (up to the above mentioned exception).

We next turn to $n\ge45$. Note that from the asymptotics above the main term from the rank is dominant. Thus main term contributes
\[
	\frac{16\sqrt3}{\sqrt{24{ (11n+d)}-1}}\sinh\left(\frac{5\pi}{66}\sqrt{24{ (11n+d)}-1}\right)\sin\pa{\pi}{11} \cos\pa{2\pi a_2}{11}.
\]
Note that $a_2$ determines whether we are in $\S_1$ or $\S_2$. Since
\[
	\cos\pa{4\pi}{11} > 0,\qquad \cos\pa{6\pi}{11} < 0,
\]
the above expresssion is positive if $a_2=2$ and negative if $a_2=3$.

Using Lemmas \ref{L:MCE} and \ref{L:rank} the error term is bounded by
\begin{multline*}
	2\CE_1{ (11n+d)}\sum_{j=1}^5 \left|\cos\pa{2\pi a_1j}{11}\right|\sin\pa{\pi j}{11} + 2\CE_2{ (11n+d)}\sum_{j=1}^5 \left|\cos\pa{2\pi a_2j}{11}\right|\sin\pa{\pi j}{11}\\
	+ 2\sum_{j=1}^5 \left|\cos\pa{2\pi a_2j}{11}\right||M_j{ (11n+d)}|.
\end{multline*}
We estimate
\begin{multline*}
	2\sum_{j=1}^5 \left|\cos\pa{2\pi a_2j}{11}\right||M_j{ (11n+d)}| \le \frac{8\sqrt3\sinh\pa{\pi\sqrt{24{ (11n+d)}-1}}{66}}{\sqrt{24{ (11n+d)}-1}}\\
	\times \max_{a_2\in\{2,3\}} \left(\frac{\sqrt{11}}{\sin\pa{\pi}{11}}\sum_{j=1}^5 \left|\cos\pa{2\pi a_2j}{11}\right|\sin\pa{\pi j}{11}+2\sin\pa{\pi}{11}\left|\cos\pa{2\pi a_2}{11}\right|\right).
\end{multline*}
The maximum is obtained for $a_2=3$ and this can be bound against 
\[
	\frac{374.1\sinh\pa{\pi\sqrt{24{ (11n+d)}-1}}{66}}{\sqrt{24{ (11n+d)}-1}}.
\]
Similarly
\[
	2\max_{a_2\in\{2,3\}} \sum_{j=1}^5 \left|\cos\pa{2\pi a_2j}{11}\right|\sin\pa{\pi j}{11}
\]
is obtained for $a_2=3$. It can be bound against $4.58$. Finally
\[
	2\max_{0\le a_1\le3} \sum_{j=1}^5 \left|\cos\pa{2\pi a_1j}{11}\right|\sin\pa{\pi j}{11}
\]
is obtained for $a_1=0$ and is bound in this case by $7$. Thus the error can be estimated against
\[
	7\CE_1{ (11n+d)} + 4.58\CE_2{ (11n+d)} + \frac{374.1\sinh\pa{\pi\sqrt{24{ (11n+d)}-1}}{66}}{\sqrt{24{ (11n+d)}-1}}.
\]
A computer calculation shows that for $a_2\in\{2,3\}$ and $n\ge45$
\begin{equation*}
	\frac{7\CE_1{ (11n+d)} + 4.58\CE_2{ (11n+d)} + \frac{374.1\sinh\pa{\pi\sqrt{24{ (11n+d)}-1}}{66}}{\sqrt{24{ (11n+d)}-1}}}{\frac{16\sqrt3}{\sqrt{24{ (11n+d)}-1}}\sinh\left(\frac{5\pi}{66}\sqrt{24{ (11n+d)}-1}\right) \sin\pa{\pi}{11}\left|\cos\pa{2\pi a_2}{11}\right|} < 1.
\end{equation*}

\subsection{$\mathcal C_{a_1,a_2}^{[5]}(n)$}

We need to show that for all $n\in\N_0$
\[
	M_{0,1}(n) \le M_{2,1}(n).
\]
This follows by Corollary 4.1 (4) of \cite{ZR}.

\end{document}